\documentclass[journal]{IEEEtran}

\usepackage{amsmath,amssymb,amsthm,mathtools,wasysym,epsfig,dsfont,color,subfigure,empheq,graphicx,tcolorbox}
\usepackage{enumerate,cite,url,algorithm, algorithmic,empheq,balance,url,multirow,etoolbox}
\usepackage{svg, multirow, optidef}
\allowdisplaybreaks

\def\eg{\emph{e.g., }}
\def\ie{\emph{i.e., }} 
\def\ea{\emph{et al. }} 
\DeclareMathOperator*{\minimize}{minimize}
\DeclareMathOperator*{\subjectto}{subject\,to}

\usepackage{tikz}
\usetikzlibrary{arrows,calc}
\usetikzlibrary{shadows}
\usetikzlibrary{decorations,decorations.pathmorphing}
\usetikzlibrary{positioning,shapes}
\tikzset{
>=stealth', % standard arrow tip
help lines/.style={dashed, thick}, % different line styles
generation lines/.style={thick},
transmission lines/.style={dotted, thick},
distribution lines/.style={densely dotted, thick},
axis/.style={<->},
important line/.style={thick},
connection/.style={thick, dotted},
}

\tikzstyle{box} = [rectangle, thick, minimum height=0.75cm, text centered, draw=black]
\tikzstyle{arrow} = [thick,->,>=stealth]
\usetikzlibrary{plotmarks}

\makeatletter
\def\thenomenclature{%
  \subsection{\nomname}
  \if@intoc\addcontentsline{toc}{subsection}{\nomname}\fi%
\nompreamble
\list{}{%
\labelwidth\nom@tempdim
\leftmargin\labelwidth
\advance\leftmargin\labelsep
\itemsep\nomitemsep
\let\makelabel\nomlabel}}
\makeatother

\title{Sensitivity-Based Vulnerability Assessment \\ of State Estimation}

\begin{document}
\author{Gonzalo~E.~Constante-Flores,~\IEEEmembership{Graduate Student~Member,~IEEE,}
		~Antonio~J.~Conejo,~\IEEEmembership{Fellow,~IEEE,}
		~Jiankang~Wang,~\IEEEmembership{Member,~IEEE}
\thanks{This work was supported in part by the National Science Foundation through grants EPCN 1808169 and ECCS 1711048.}
\thanks{G. E. Constante-Flores and J.K. Wang are with the Department of Electrical and Computer Engineering, The Ohio State University, Columbus, OH, 43210 USA (e-mail: \{constanteflores.1,wang.6536\}@osu.edu).}
\thanks{A. J. Conejo is with the Department of Integrated Systems Engineering and with the Department of Electrical and Computer Engineering, The Ohio State University, Columbus, OH, 43210 USA (e-mail: conejo.1@osu.edu).}}

\maketitle

\begin{abstract}
We propose a technique to assess the vulnerability of the power system state estimator. We aim at identifying measurements that have a high potential of being the target of false data injection attacks. From an adversary's point of view, such measurements have to show the following characteristics: i) being influential on the variable estimates, and ii) corrupting their measured values is likely to be undetected. Additionally, such characteristics should not change significantly with the system's operating condition. Our technique provides a systematic way of identifying measurements with such characteristics. We illustrate our methodology on a 4-bus system and the New England 39-bus system.
\end{abstract}
\section*{Nomenclature}
\addcontentsline{toc}{section}{Nomenclature}
\hspace{-1em}\textit{Sets}
\begin{IEEEdescription}[\IEEEusemathlabelsep\IEEEsetlabelwidth{$r^\textrm{up}R^\textrm{down}$}\setlength{\IEEElabelindent}{0pt}]
\item[$\mathcal{N}$] Set of buses indexed by $i$.
\item[$\mathcal{N}_i$] Set of buses connected to bus $i$.
\item[$\mathcal{V}$] Set of buses with voltage magnitude measurements.
\item[$\mathcal{P}/\mathcal{Q}$] Set of buses with active/reactive power measurements.
\item[$\mathcal{Z}$] Set of buses with zero injections.
\item[$\mathcal{P}_f/\mathcal{Q}_f$] Set of branches with active/reactive power flow measurements.
\end{IEEEdescription}
\textit{Parameters}
\begin{IEEEdescription}[\IEEEusemathlabelsep\IEEEsetlabelwidth{$r^\textrm{up}R^\textrm{down}$}\setlength{\IEEElabelindent}{0pt}]
\item[$w_i^x$] Weighting factor for a measurement at bus $i$, where superindex $x=V,P,Q,$ refers to voltage, active power, and reactive power, respectively.
\item[$w_{ij}^x$] Weighting factor for a measurement at line $ij$, where superindex $x=P,Q,$ refers to active power and reactive power flows, respectively.
\item[$v_i^m$] Voltage magnitude measurement at bus $i$.
\item[$P_i^m/Q_i^m$] Active/reactive power injection measurement at bus $i$.
\item[$P_{ij}^m/Q_{ij}^m$] Active/reactive power flow measurement from bus $i$ to $j$.
\item[$G_{ij}/B_{ij}$] Real/imaginary part of the $ij$th entry of the admittance matrix.
\item[$b_{ij}^\text{sh}$] Shunt susceptance of line $ij$.
\end{IEEEdescription}
\textit{Variables}
\begin{IEEEdescription}[\IEEEusemathlabelsep\IEEEsetlabelwidth{$r^\textrm{up}R^\textrm{down}$}\setlength{\IEEElabelindent}{0pt}]
\item[$v_i/\theta_i$] Voltage magnitude/angle at bus $i$.
\item[$P_i / Q_i$] Active/reactive power injection at bus $i$.
\item[$P_{ij} / Q_{ij}$] Active/reactive power flow $ij$.
\end{IEEEdescription}

\section{Introduction}

\subsection{Motivation and aim}
%Energy Management Systems (EMS) is a decision-support system, which helps the system's operator to control and optimize the operation of the power grid, and it comprises four main functionalities: network model building, security assessment, automatic generation control, and economic dispatch. Network model building function aims at determining the topology of the system and estimating the operating condition of the system, which is known as state estimation.

One of the key functions of Energy Management Systems (EMS) is state estimation (SE), which aims at finding the most likely estimate of the system state (\ie voltage phasors) given the network topology and parameters, and a set of real-time measurements from remote meters and meters \cite{Abur2004}. Although such measurements always contain small errors due to the accuracy of the corresponding meters, they can also contain gross errors due to failures in telemetry and/or meters. Gross errors can also be intentionally injected to deceive the control decisions of the EMS functions (\ie security assessment, automatic generation control, and economic dispatch) by exploiting the vulnerabilities of the telemetry systems to cyber-attacks.

Such cyber-attacks, known as data-integrity attacks, aim at altering the breaker status and/or measurements while remaining undetected \cite{Deng2017,Liang2017}. In particular, we focus on a class of attacks known as false data injection attacks (FDIAs) where an adversary compromises a small subset of analog measurements (\eg voltage magnitude and angles, power flows, and power injections) to conceal a particular goal, \eg deceive system operator to change generation dispatch, congest transmission lines, or produce cascading failures \cite{Liu2009}. 

Within the context above, we aim at analyzing the vulnerabilities of the state estimation against FDIAs based on sensitivity analysis. Such vulnerabilities are characterized in terms of the chance of an attack to significantly influence (if perturbed) the optimal estimates while remaining undetected. 

\subsection{Related work}
The theoretical framework for sensitivity analysis in nonlinear optimization used in this paper is stated in \cite{Castillo2006a} and \cite{Castillo2006}, following the pioneering works of Fiacco \cite{Fiacco1983}, Enevoldsen \cite{Enevoldsen1994}, and Bonnans and Shapiro \cite{Bonnans2013}. Power systems applications of this sensitivity analysis framework include state estimation \cite{Minguez2007} and \cite{Caro2011} and pricing \cite{Conejo2005}. Other engineering applications include reliability analysis \cite{Castillo2003}, calculus of variations \cite{Castillo2008} and civil infrastructure optimal design \cite{Minguez2006}. 

Given the characteristics of the nonlinear state estimator \cite{Liang2017}, its vulnerability has been studied only in a few works. The vulnerability of the state estimator has been quantified by the minimum number of sensors that have to be compromised to stage a stealthy FDIA, which can be formulated as a minimum cardinality problem \cite{Kosut2010,Dan2010}. Hug and Giampapa \cite{Hug2012} propose a graph-based algorithm to find the set of compromised sensors needed to stage an unobservable attack assuming that the adversary has perfect information of the system. Rahman \ea \cite{Rahman2013} extend the work in \cite{Hug2012} by considering incomplete information of the system. Zhao \ea \cite{Zhao2018a} propose a framework to analyze the vulnerability of the nonlinear state estimation from the system operator's point of view and propose countermeasures. Jin \ea \cite{Jin2019} formulate a framework based on a semi-definite convexification of the FDIA to find a near-optimal attack strategy and analyze the attack stealthiness. They provide theoretical guarantees of sparsity and unobservability. However, this formulation depends on the adversaries' objective, which is not necessarily always available to the system's operator.

\subsection{Contribution and paper organization}
In this paper, we tailor the sensitivity analysis methodology in \cite{Castillo2006} to efficiently analyze vulnerabilities in the state estimation problem with respect to FDIAs. Unlike the existing literature, our methodology does not depend on the adversary's objective nor quantify the vulnerability by the minimum number of sensors needed to be compromised to stage an unobservable FDIA. From the system's operator point of view, we rather focus on the vulnerability based on endogenous factors of the state estimation and the power grid, such as the measurement configuration, system topology, and network parameters. The contributions of this work are threefold:

\begin{itemize}
    \item To analyze the vulnerabilities in the state estimation problem with respect to FDIAs based on sensitivity analysis. We identify such vulnerabilities in terms of the stealthiness and impactfulness characteristics of a FDIA when it targets a particular measurement. The sensitivity analysis methodology allows us to compute both characteristics of all the measurements simultaneously.
    \item To propose three scores to quantify and rank the vulnerability of each measurement to FDIAs. These scores allow us to quantify the vulnerabilities of the measurements, which can help to identify vulnerable areas of the system and to improve its security.
    \item To assess the variations of the sensitivities with respect to different operating conditions based on a singular value decomposition (SVD) approach. We aim at identifying whether the vulnerabilities of the state estimation vary with respect to the operating condition of the system or they remain almost invariant. The latter case would imply that the vulnerabilities are mainly dependent on the network topology and its parameters, and the configuration of the measurements.
\end{itemize}

Although we illustrate our methodology in the weighted least squared (WLS) state estimator, such methodology can be implemented using other estimators (\eg robust estimators) as long as they can be stated as a continuous optimization problem and their solution holds the Karush-Kuhn-Tucker (KKT) optimality conditions.

The remainder of this paper is organized as follows. In
Section II, we present the characterization of vulnerable measurements, state the estimation formulation and the analytical expressions to compute the sensitivities. The method to identify if such sensitivities change with the system's operating conditions is described in Section III. The proposed model is validated through numerical experiments in two test systems in Section IV. The main conclusions of the paper are summarized in Sections V.

\section{Vulnerability Analysis}
In this section, we characterize the vulnerability of the measurements against FDIA. Also, we present the state estimation formulation and derive the analytical expressions to compute the sensitivities of the objective and estimated variables with respect to parameters and measurements.

\subsection{Characterizing vulnerable measurements}

The goal of a FDIA is to stealthily modify measurements to introduce gross errors in the variable estimates, which are then used in other control applications (\eg security-constrained optimal power flow and security analysis) \cite{Liang2017}. This goal shows two main characteristics:  

\subsubsection{Stealthiness}
Once the solution of the state estimator is computed, gross errors are detected by comparing the sum of squared errors with a bad data detection (BDD) flag. In the case of the WLS estimator, the widely adopted criteria for this flag comes from a $\chi^2$ distribution \cite{Abur2004,Gomez2018}. Note that if the state estimator is formulated as an optimization problem, the sum of squared errors is the value of the objective function. 
    
An adversary aims at modifying measurements without triggering the BDD flag, which could hinder the successful staging of the attack. Thus, an attacker would look to corrupt the measurements that do not change significantly the objective function when they are perturbed, which means that the rate of change of the objective function with respect to the measurement is small.

\subsubsection{Impactfulness}
Besides remaining undetected, an adversary aims at causing a large change in the variable estimates without significantly modifying the measurement under attack \ie the rate of change in the variable estimate as a measurement changes has to be large. Since the state estimation can be also understood as a nonlinear regression problem, this characteristic turns out to be the definition of leverage point in regression analysis \cite{Laurent1992}. 
Measurements with high leverage have two important characteristics: i) they have a significant influence on the variable estimate when they are perturbed, and ii) they can be eliminated without losing system observability unless they are critical measurements \cite{Abur2004}.

A measurement showing both characteristics is a high-potential target for cyber-attack as an adversary can stage an impactful attack while remaining likely undetected. We note that both characteristics can be described in terms of the sensitivities of the objective function and the variable estimates with respect to the measurements. The remainder of this section presents a technique to systematically compute both sensitivities for all the measurements simultaneously using solely state estimation information.

\subsection{State Estimation Formulation}
The weighted least squares (WLS) state estimation can be formulated as an equality-constrained optimization problem as follows:

{\small{
\vspace{-0.2cm}
\begin{subequations}
\begin{flalign}
            &\minimize_{\boldsymbol{x}}\quad 
            \smashoperator{\sum_{i\in \mathcal{V}}} w_i^V (v_i^m-v_i)^2 \nonumber \\
            &\phantom{\minimize_{\boldsymbol{x}}\quad } + \smashoperator{\sum_{i\in \mathcal{P}}} w_i^P (P_i^m-P_i)^2 +\smashoperator{\sum_{i\in \mathcal{Q}}} w_i^Q (Q_i^m-Q_i)^2 \nonumber \\
            &\phantom{\minimize_{\boldsymbol{x}}\quad } + \smashoperator{\sum_{(i,j)\in \mathcal{P}_f}} w_{ij}^P (P_{ij}^m-P_{ij})^2 + \smashoperator{\sum_{(i,j)\in \mathcal{Q}_f}} w_{ij}^Q (Q_{ij}^m-Q_{ij})^2 \label{eq:1a} \\
            & \subjectto \nonumber \\
            & P_i = v_i \smashoperator{\sum_{j\in\mathcal{N}_i}}v_j(G_{ij}\text{cos}(\theta_{ij})+B_{ij}\text{sin}(\theta_{ij})),\,i \in \mathcal{P}    \hspace{-0.5cm} \label{eq:1b} \\
            & Q_i = v_i \smashoperator{\sum_{j\in\mathcal{N}_i}}v_j(G_{ij}\text{sin}(\theta_{ij})-B_{ij}\text{cos}(\theta_{ij})),\,i \in \mathcal{Q}  \hspace{-0.5cm} \label{eq:1c} \\
            & P_{ij} = v_iv_j(G_{ij}\text{cos}(\theta_{ij})+B_{ij}\text{sin}(\theta_{ij}))-G_{ij}v_i^2,\,(i,j) \in \mathcal{P}_f   \hspace{-0.5cm} \label{eq:1d} \\
            & Q_{ij} = v_iv_j(G_{ij}\text{sin}(\theta_{ij})-B_{ij}\text{cos}(\theta_{ij})) \nonumber\\
            & \hspace{3.5cm} +v_i^2(B_{ij}-b_{ij}^\textrm{sh}/2),\,(i,j) \in \mathcal{Q}_f \hspace{-0.5cm} \label{eq:1e} \\
            & 0 = v_i \smashoperator{\sum_{j\in\mathcal{N}_i}}v_j(G_{ij}\text{cos}(\theta_{ij})+B_{ij}\text{sin}(\theta_{ij})),\,i \in \mathcal{Z}    \hspace{-0.5cm} \label{eq:1f} \\
            & 0 = v_i \smashoperator{\sum_{j\in\mathcal{N}_i}}v_j(G_{ij}\text{sin}(\theta_{ij})-B_{ij}\text{cos}(\theta_{ij})),\,i \in \mathcal{Z}.    \hspace{-0.5cm} \label{eq:1g}
\end{flalign}
\end{subequations}}}

The objective \eqref{eq:1a} is to minimize the weighted sum of squared errors. Constraints \eqref{eq:1b}-\eqref{eq:1c} represent the active and reactive power injections of the buses with available injection measurements. Constraints \eqref{eq:1d}-\eqref{eq:1e} represent the active and reactive power flows of the lines with available flow measurements. Constraints \eqref{eq:1f}-\eqref{eq:1g} correspond to the zero-injections (\ie exact pseudo-measurements).

The above problem can be expressed in compact form as:
\begin{subequations}
\begin{align}
\minimize_{\boldsymbol{x}}\;\;&J(\boldsymbol{x},\boldsymbol{a},\boldsymbol{z}) \label{eq:2a} \\
\text{subject to}\;\;&\boldsymbol{c}(\boldsymbol{x},\boldsymbol{a})=0: \quad \boldsymbol{\lambda}, \label{eq:2b}
\end{align} \label{eq:2}
\end{subequations}
where $J$ is the sum of squared errors, $\boldsymbol{x}\in \mathbb{R}^n$ denotes all optimization variables $(v_i,\theta_i,P_i,Q_i,P_{ij},Q_{ij})$, $\boldsymbol{z}\in \mathbb{R}^p$ represents the measurements $(v_i^m,\theta_i^m,P_i^m,Q_i^m,P_{ij}^m,Q_{ij}^m)$, $\boldsymbol{a}\in \mathbb{R}^q$ represents all parameters ($w_i^V,w_i^P,w_i^Q,w_{ij}^P,w_{ij}^Q, G_{ij},B_{ij},b_{ij}^\text{sh}$),  and $\boldsymbol{\lambda}\in \mathbb{R}^r$ denotes the Lagrange multiplier vector. Note that the equality constraints only depend on the optimization variables and the parameters, not the measurements. 

\subsection{Feasible Perturbations and Sensitivity Analysis}
Let $\boldsymbol{x}^\ast$ be a local optimal solution of Problem \eqref{eq:2}, and assume that $\boldsymbol{x}^\ast$ is regular \ie the constraint gradients $\nabla c_k(\boldsymbol{\boldsymbol{x}^\ast}),\,k=1,\dots,r$, are linearly independent \cite{Luenberger2008}. Then, the KKT first-order optimality conditions are \cite{Luenberger2008}
\begin{align}
    & \nabla_{\boldsymbol{x}}J(\boldsymbol{x}^\ast,\boldsymbol{a},\boldsymbol{z}) + \sum_{k=1}^{r}\lambda_k^\ast \nabla_{\boldsymbol{x}}  c_k(\boldsymbol{x}^\ast,\boldsymbol{a}) = 0 \label{eq:3}\\
    & c_k(\boldsymbol{x}^\ast,\boldsymbol{a})=0,\; k = 1,\dots,r, \label{eq:4}
\end{align}
where conditions \eqref{eq:4} are the primal feasibility conditions.

To determine the sensitivity equations with respect to the parameters and measurements, we perturb $\boldsymbol{x}^\ast,\boldsymbol{\lambda}^\ast,J^\ast,\boldsymbol{a},\boldsymbol{z}$ in such a way that the KKT conditions still hold \cite{Castillo2006}. Thus, to obtain such equations, we differenciate the objective function \eqref{eq:2a} and the optimality conditions \eqref{eq:3}-\eqref{eq:4} as follows:
\begin{gather}
     \left[\nabla_{\boldsymbol{x}}J(\boldsymbol{x}^\ast,\boldsymbol{a},\boldsymbol{z})\right]^\top d\boldsymbol{x} + 
     \left[\nabla_{\boldsymbol{a}}J(\boldsymbol{x}^\ast,\boldsymbol{a},\boldsymbol{z})\right]^\top d\boldsymbol{a}+ \nonumber \\
     \left[\nabla_{\boldsymbol{z}}J(\boldsymbol{x}^\ast,\boldsymbol{a},\boldsymbol{z})\right]^\top d\boldsymbol{z}  -dJ = 0 \label{eq:5}\\
     \left[\nabla_{\boldsymbol{xx}}J(\boldsymbol{x}^\ast,\boldsymbol{a},\boldsymbol{z}) + \sum_{k=1}^{r}\lambda_k^\ast \nabla_{\boldsymbol{xx}} c_k(\boldsymbol{x}^\ast,\boldsymbol{a})\right] d\boldsymbol{x} + \nonumber \\  
     \left[\nabla_{\boldsymbol{xa}}J(\boldsymbol{x}^\ast,\boldsymbol{a},\boldsymbol{z}) + \sum_{k=1}^{r}\lambda_k^\ast \nabla_{\boldsymbol{xa}} c_k(\boldsymbol{x}^\ast,\boldsymbol{a})\right] d\boldsymbol{a} + \nonumber \\
     \nabla_{\boldsymbol{xz}}J(\boldsymbol{x}^\ast,\boldsymbol{a},\boldsymbol{z}) d\boldsymbol{z} + \nabla_{\boldsymbol{x}} \boldsymbol{c}(\boldsymbol{x}^\ast,\boldsymbol{a}) d \boldsymbol{\lambda} = 0 \label{eq:6} \\
     \left[\nabla_{\boldsymbol{x}} \boldsymbol{c}(\boldsymbol{x}^\ast,\boldsymbol{a})\right]^\top d\boldsymbol{x} + \left[\nabla_{\boldsymbol{a}} \boldsymbol{c}(\boldsymbol{x}^\ast,\boldsymbol{a})\right]^\top d\boldsymbol{a}= 0 \label{eq:7}
\end{gather}

The above system of equations can be expressed in matrix form as follows:
\begin{gather}
\begin{bmatrix}
\boldsymbol{J}_{\boldsymbol{x}} & \boldsymbol{J}_{\boldsymbol{a}} & \boldsymbol{J}_{\boldsymbol{z}}  & \boldsymbol{0} & -1 \\
\boldsymbol{J}_{\boldsymbol{xx}} & \boldsymbol{J}_{\boldsymbol{xa}} & \boldsymbol{J}_{\boldsymbol{xz}} & \boldsymbol{C}_{\boldsymbol{x}}^\top & 0 \\
\boldsymbol{C}_{\boldsymbol{x}} & \boldsymbol{C}_{\boldsymbol{a}} & \boldsymbol{0}& \boldsymbol{0} & 0 \\
\end{bmatrix}  
\begin{bmatrix}
d\boldsymbol{x} \\
d\boldsymbol{a} \\
d\boldsymbol{z} \\
d\boldsymbol{\lambda} \\
dJ
\end{bmatrix} = \boldsymbol{0} \label{eq:8}
\end{gather}
where the vectors and submatrices are defined in Appendix A.

Then, equation \eqref{eq:5} can be written as:
\begin{gather}
    \boldsymbol{T} \begin{bmatrix} d\boldsymbol{x} & d\boldsymbol{\lambda} & dJ \end{bmatrix}^\top = \boldsymbol{S}_{\boldsymbol{a}} d\boldsymbol{a}+\boldsymbol{S}_{\boldsymbol{z}} d \boldsymbol{z}, \label{eq:9}
\end{gather}
 where the matrices $\boldsymbol{T}$, $\boldsymbol{S}_{\boldsymbol{a}}$, and $\boldsymbol{S}_{\boldsymbol{z}}$ are
\begin{align*}
    \boldsymbol{T} =  
    \begin{bmatrix}
        \boldsymbol{J}_{\boldsymbol{x}} & \boldsymbol{0} & -1 \\
        \boldsymbol{J}_{\boldsymbol{xx}} & \boldsymbol{C}_{\boldsymbol{x}}^\top & \boldsymbol{0} \\
        \boldsymbol{C}_{\boldsymbol{x}} & \boldsymbol{0} & \boldsymbol{0}
    \end{bmatrix}, \\
    \boldsymbol{S}_{\boldsymbol{a}}^\top =  -\begin{bmatrix} \boldsymbol{J}_{\boldsymbol{a}} & \boldsymbol{J}_{\boldsymbol{xa}} & \boldsymbol{J}_{\boldsymbol{a}} \end{bmatrix}, \\
    \boldsymbol{S}_{\boldsymbol{z}}^\top =  -\begin{bmatrix} \boldsymbol{J}_{\boldsymbol{z}} & \boldsymbol{J}_{\boldsymbol{xz}} & \boldsymbol{J}_{\boldsymbol{z}} \end{bmatrix}.
\end{align*}
Therefore, we can express \eqref{eq:9} as follows:
\begin{gather}
\begin{bmatrix} d\boldsymbol{x} & d\boldsymbol{\lambda} & dJ \end{bmatrix}^\top = \boldsymbol{T}^{-1}\boldsymbol{S}_{\boldsymbol{a}}d\boldsymbol{a} + \boldsymbol{T}^{-1}\boldsymbol{S}_{\boldsymbol{z}}d\boldsymbol{z}, \label{eq:10}
\end{gather}
which can be solved using the superposition principle by replacing $d\boldsymbol{z}$ and $d\boldsymbol{a}$ by the $p$-dimensional and $q$-dimensional identity matrices, respectively. Then, we obtain the matrices with all sensitivities with respect to parameters and measurements
\begin{gather}
\begin{bmatrix} \dfrac{\partial \boldsymbol{x}}{\partial \boldsymbol{a}} & \dfrac{\partial \boldsymbol{\lambda}}{\partial \boldsymbol{a}}\end{bmatrix}^\top = -\boldsymbol{H}_{\boldsymbol{x}}^{-1}\boldsymbol{H}_{\boldsymbol{a}}, \quad \frac{\partial J}{\partial \boldsymbol{a}} = \boldsymbol{J}_{\boldsymbol{a}}+ \boldsymbol{J}_{\boldsymbol{x}}\frac{\partial \boldsymbol{x}}{\partial \boldsymbol{a}} \label{eq:11}\\
\begin{bmatrix} \dfrac{\partial \boldsymbol{x}}{\partial \boldsymbol{z}} & \dfrac{\partial \boldsymbol{\lambda}}{\partial \boldsymbol{z}}\end{bmatrix}^\top = -\boldsymbol{H}_{\boldsymbol{x}}^{-1}\boldsymbol{H}_{\boldsymbol{z}}, \quad \frac{\partial J}{\partial \boldsymbol{z}} = \boldsymbol{J}_{\boldsymbol{z}}+ \boldsymbol{J}_{\boldsymbol{x}}\frac{\partial \boldsymbol{x}}{\partial \boldsymbol{z}}, \label{eq:12}
\end{gather}
where
\begin{align*}
    \boldsymbol{H}_{\boldsymbol{x}} =  
    \begin{bmatrix}
        \boldsymbol{J}_{\boldsymbol{xx}} & \boldsymbol{C}_{\boldsymbol{x}}^\top \\
        \boldsymbol{C}_{\boldsymbol{x}} & \boldsymbol{0}
    \end{bmatrix}, \quad
    \boldsymbol{H}_{\boldsymbol{a}} =  \begin{bmatrix} \boldsymbol{J}_{\boldsymbol{xa}} \\ \boldsymbol{C}_{\boldsymbol{a}} \end{bmatrix}, \quad
    \boldsymbol{H}_{\boldsymbol{z}} =  \begin{bmatrix} \boldsymbol{J}_{\boldsymbol{xz}} \\ \boldsymbol{0} \end{bmatrix}.
\end{align*}

Clearly, the sensitivities of the objective and the variable estimates with respect to the measurements, which allow us to define the vulnerability of each measurement, can be computed from \eqref{eq:12}. 

\subsection{Identifying vulnerable measurements}

To better visualize the stealthiness and impactfulness of a measurements $z_{\ell}$, we propose three scores to rank the vulnerability of a measurement $z_{\ell}$: i) $\text{S-score}(z_{\ell})$, which quantifies how likely is a FDIA to be undetected, ii) $\text{L-score}(z_{\ell})$, which quantifies the influence of a FDIA in the variables estimates, and iii) $\text{V-score}(z_{\ell})$, which is a convex combination of the previous scores. These scores are defined as follows: 
\begin{gather}
    \text{S-score}(z_{\ell}):=
    f\left(\gamma^{-\frac{\left\lvert z_{\ell}\frac{\partial J}{\partial z_{\ell}}\right\rvert}{{\max\limits_{1\le k\le p} \left\{ \left\lvert z_k \frac{\partial J}{\partial z_k}\right\rvert\right\} }}}\right), \label{eq:13} \\
\text{L-score}(z_{\ell}):= g\left(\left\lVert\frac{\partial \boldsymbol{x}}{\partial z_{\ell}}\right\rVert \right), \label{eq:14} \\
    \text{V-score}(z_{\ell}):= \alpha\left(\text{S-score}(z_{\ell})\right) + (1-\alpha)\left(\text{L-score}(z_{\ell})\right), \label{eq:15}
\end{gather}
where $\gamma >0$, $\alpha \in [0,1]$, and $f(\cdot)$ and $g(\cdot)$ are non-decreasing functions with range and domain on $\left[0,1\right]$. It is relevant to note that to score the leverage of a measurement $z_\ell$, we consider the norm of the sensitivities of the variable estimates with respect to it $\left\lVert\frac{\partial \boldsymbol{x}}{\partial z_{\ell}}\right\rVert$. This allows us to take into account the influence of such measurement not only on its  corresponding variable estimate (\ie self-sensitivity), but also its impact on the other variable estimates. 

\begin{figure}[t]
    \centering
    \includegraphics[width = 0.48\textwidth]{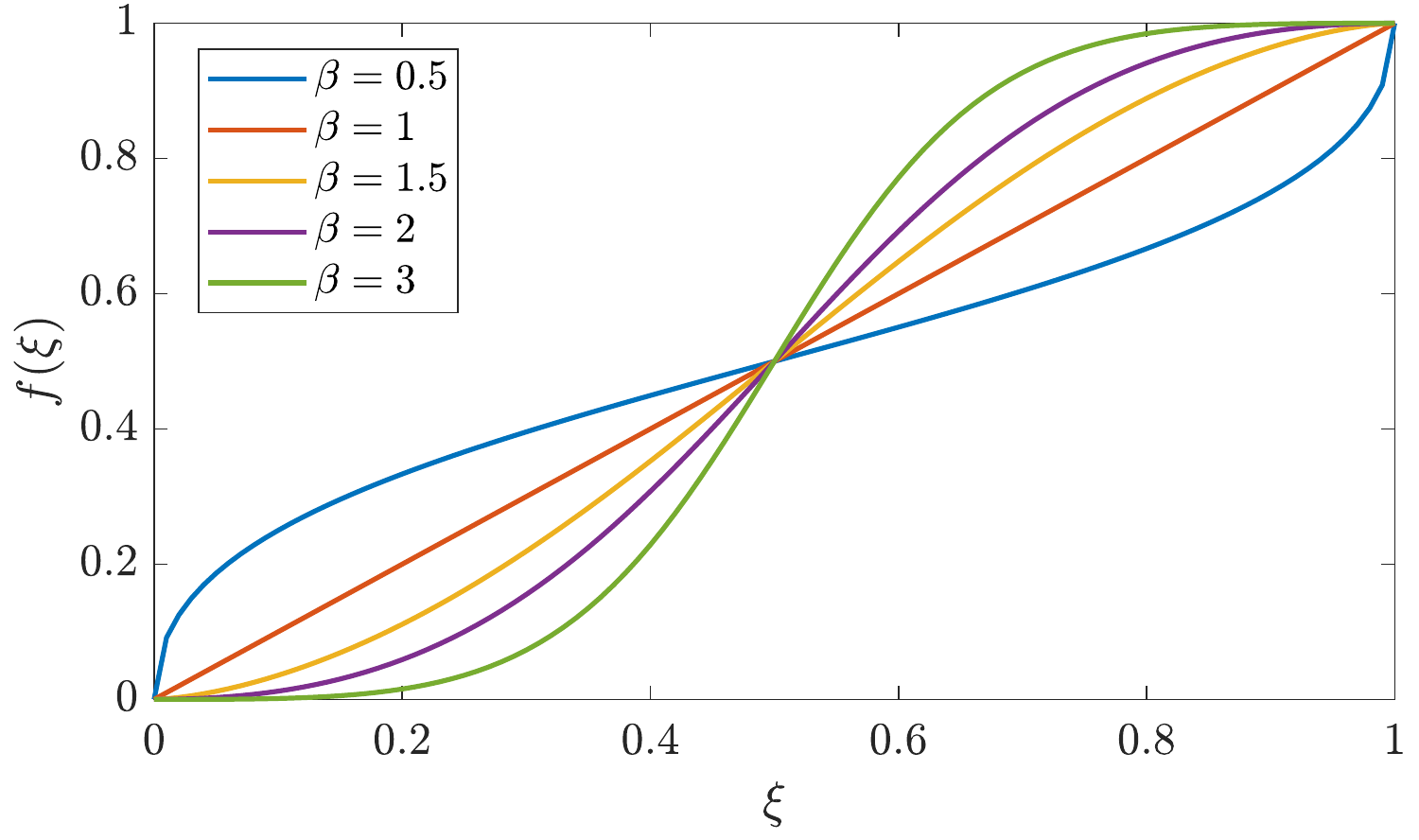}
    \caption{S-shaped function.}
    \label{fig:S_function}
\end{figure}

The proposed scores are closer to 1 when a measurement is more vulnerable. It is noteworthy that $f(\cdot)$ and $g(\cdot)$ and their arguments are user-defined. We suggest an S-shaped function for both scores, such as
\begin{align*}
f(\xi) = \begin{cases}
0, & \xi \le 0 \\
\frac{1}{1+ \left(\frac{\xi}{1-\xi} \right)^{-\beta}}, & 0\le \xi \le 1 \\
1, & \xi \ge 1 \\
\end{cases}
\end{align*}
where $\beta>0$. Note that smaller values of $\beta$ render more conservative scores as the function rapidly downweight the scores as they distance from 1 as depicted in Fig.~\ref{fig:S_function}.

Finally, we note that the procedure to compute the sensitivities with respect to measurements and the proposed scores is summarized in Algorithm~\ref{alg:SA_SE}.

\begin{algorithm}[ht]
\caption{Sensitivity Analysis of State Estimation}
\label{alg:SA_SE}
\begin{algorithmic}[1]
\REQUIRE Optimal solution $\left(J^\ast, \boldsymbol{x}^\ast, \boldsymbol{\lambda}^\ast \right)$, parameters $\beta > 0$, $\gamma >0$, and $\alpha \in [0,1]$, functions $f(\cdot)$ and $g(\cdot)$.
\STATE \textit{Compute.} $\boldsymbol{H}_{\boldsymbol{x}},\boldsymbol{H}_{\boldsymbol{z}},\boldsymbol{J}_{\boldsymbol{x}},\boldsymbol{J}_{\boldsymbol{z}}$.
\STATE \textit{Compute.} $\frac{\partial J}{\partial \boldsymbol{z}}$ and $\frac{\partial \boldsymbol{x}}{\partial \boldsymbol{z}}$ from \eqref{eq:12}.
\FOR{$\ell = 1,\dots,p$}
\STATE $\text{S-score}(z_{\ell}) \gets \text{Evaluate } \eqref{eq:13}$.
\STATE $\text{L-score}(z_{\ell}) \gets \text{Evaluate } \eqref{eq:14}$.
\STATE $\text{V-score}(z_{\ell}) \gets \text{Evaluate } \eqref{eq:15}$.
\ENDFOR
\ENSURE Sensitivities $\left(\frac{\partial J}{\partial z_{\ell}},\frac{\partial \boldsymbol{x}}{\partial z_{\ell}}\right)$, and Scores $\big(\text{S-score}(z_{\ell}),$ $\text{L-score}(z_{\ell}),$ $\text{V-score}(z_{\ell})\big),\, \ell = 1,\dots,p$.
\end{algorithmic}
\end{algorithm}

%Our suggested scores are based on the idea of universal threshold in the lasso variable selection \cite{Gao2018}. In the case of $\text{S-score}(z_{\ell})$ (resp. $\text{L-score}(z_{\ell})$), the constant $\textrm{log}(p)$ in both scores can be chosen to be smaller (resp. larger) to be more conservative.

%To identify the measurements with potential to be targets of cyber-attacks, we define the following index sets:
%\begin{gather}
%    \mathcal{D} = \left\{\ell\,:\,\text{S-score}(z_{\ell}) \ge 0\right\}, \label{eq:16} \\
%    \mathcal{L} = \left\{\ell\,:\,\text{L-score}(z_{\ell}) \ge 0\right\}, \label{eq:17} \\
%    \mathcal{T} = \left\{\ell\,:\, \ell \in \mathcal{D} \cap \mathcal{L}\right\}, \label{eq:18}
%\end{gather}
%where $\mathcal{D}$ is the set of measurements with likely undetectable FDIA, $\mathcal{L}$ denotes the set of measurements that are highly influential in their variable estimates, and $\mathcal{T}$ is the set of measurements with high-potential of being target for cyber-attacks.

\section{Robustness Analysis}
In this section, we present a method to identify whether or no the sensitivities change with the system's operating condition. 
\subsection{Preprocessing}
To determine if the sensitivity vectors show significant changes with respect to the operating points, we consider $t$ different operating conditions and compute their corresponding sensitivities. Then, we arrange these sensitivities in matrices $\boldsymbol{X} \in \mathbb{R}^{t \times np}$ and $\boldsymbol{J} \in \mathbb{R}^{t \times p}$  as follows:
\begin{align}
    \boldsymbol{X} = \begin{bmatrix} $-----$\;\boldsymbol{x}_1^\top\;$-----$  \\ $-----$\;\boldsymbol{x}_2^\top\;$-----$   \\ \vdots  \\ $-----$\;\boldsymbol{x}_t^\top\;$-----$    \end{bmatrix}, \; \boldsymbol{J} = \begin{bmatrix} $-----$\;\boldsymbol{J}_1^\top\;$-----$  \\ $-----$\;\boldsymbol{J}_2^\top\;$-----$   \\ \vdots  \\ $-----$\;\boldsymbol{J}_t^\top\;$-----$  \end{bmatrix}, \label{eq:19}
\end{align}
where $ \boldsymbol{x}_k:= \text{vec}\left( \frac{\partial \boldsymbol{x}}{\partial \boldsymbol{z}} \right)$ and $ \boldsymbol{J}_k:= \frac{\partial J}{\partial \boldsymbol{z}}$ are the sensitivities at a given operating condition $k$. Each column of $\boldsymbol{X}$ and $\boldsymbol{J}$ correspond to a particular sensitivity for all the operating conditions.

Note that in \eqref{eq:19} we assume that every sensitivity vector $\boldsymbol{x}_k$ and $\boldsymbol{J}_k$ has the same dimension \ie the system topology and measurement configuration remain unaltered, which might not be always true. If the dimensions of the sensitivity vectors are different, then it is necessary to keep only the sensitivities that are common for all the operating conditions.

Singular value decomposition (SVD) allow us to determine if such sensitivities significantly vary depending on the different operating points. Before computing the SVD of both matrices, it is necessary to subtract the mean of each column (\ie each column has zero-mean). We compute the row vectors containing the means of every column as
\begin{align}
    \bar{\boldsymbol{x}} = \dfrac{1}{t} \sum\limits_{k=1}^t \boldsymbol{x}_k,\quad \bar{\boldsymbol{j}} = \dfrac{1}{t} \sum\limits_{k=1}^t \boldsymbol{J}_k, \label{eq:20}
\end{align}
where $\bar{\boldsymbol{x}} \in \mathbb{R}^{1\times np}$ and $\bar{\boldsymbol{j}} \in \mathbb{R}^{1 \times p}$. 
Then, we can compute the elements of the centered matrices as follows:
\begin{align}
    \tilde{\boldsymbol{X}} = \boldsymbol{X}-\boldsymbol{1}_{\boldsymbol{p}}\cdot\bar{\boldsymbol{x}},\quad
    \tilde{\boldsymbol{J}} = \boldsymbol{J}-\boldsymbol{1}_{\boldsymbol{p}}\cdot\bar{\boldsymbol{j}}, \label{eq:21}
\end{align}
where $\boldsymbol{1}_{\boldsymbol{p}}$ is the $t$-dimensional all-ones column vector, and $\tilde{\boldsymbol{X}}$ and $\tilde{\boldsymbol{J}}$ are the mean-centered sensitivity matrices.

\subsection{Singular Value Decomposition}
The SVD allows us to determine if the sensitivities are significantly affected by the different operating points. We compute the SVD of both standardized matrices as follows:
\begin{align}
    \tilde{\boldsymbol{X}} = \boldsymbol{U} \boldsymbol{\Sigma} \boldsymbol{V}^\top = \sum_{i=1}^{\min \{t,np\}} \sigma_i \boldsymbol{u}_i \boldsymbol{v}^\top_i, \label{eq:22}\\ 
    \tilde{\boldsymbol{J}} = \hat{\boldsymbol{U}} \hat{\boldsymbol{\Sigma}} \hat{\boldsymbol{V}}^\top = \sum_{i=1}^{\min \{t,p\}} \hat{\sigma}_i \hat{\boldsymbol{u}}_i \hat{\boldsymbol{v}}_i^\top,  \label{eq:23}
\end{align}
where the diagonal elements of $\boldsymbol{\Sigma}$ and $\hat{\boldsymbol{\Sigma}}$ are the singular values of $\tilde{\boldsymbol{X}}$ and $\tilde{\boldsymbol{J}}$, respectively, and they are ordered from largest to smallest. 

If the largest singular values are significantly larger than the smallest ones, we can say that the sensitivities are not strongly dependent on the system's operating condition.
Such characteristic is key for a cyber-attack because it means that the sensitivities depend on factors that do not change significantly over time (\eg system's topology, lines' parameters, and measurements' locations and precisions). Thus, it allows the adversary to identify the target measurements off-line, and to stage an attack on one of these measurements without knowing other measurements.

To quantify the proportion of the variance of the mean-centered sensitivity matrices $\tilde{\boldsymbol{X}}$ and $\tilde{\boldsymbol{J}}$ that is captured by their first $r$ singular values, we define the cumulative energy (CE) as follows:
\begin{align}
    \text{CE}(\tilde{\boldsymbol{X}};r) := \frac{\hspace{0.3cm} \smashoperator{\sum\limits_{i=1}^r} \sigma_i \hspace{0.3cm}}{\hspace{0.3cm} \smashoperator{\sum\limits_{i=1}^{\min \{t,np\}}} \sigma_i}, \;
    \text{CE}(\tilde{\boldsymbol{J}};r) := \frac{\hspace{0.3cm} \smashoperator{\sum\limits_{i=1}^r} \hat{\sigma}_i \hspace{0.3cm}}{\hspace{0.3cm}\smashoperator{\sum\limits_{i=1}^{\min \{t,p\}}} \hat{\sigma}_i} \label{eq:24}
\end{align}

\section{Illustrative Examples}
In this section, two case studies are analyzed considering a 4-bus system and the New England 39-bus system. %The data for the 4-bus system are presented in Appendix B whereas the data for the New England 39-bus system can be retrieved from MATPOWER \cite{Zimmerman2011}. 
The weights of the voltage measurements are assumed to be $w^V = 1\text{e}^{4}$ whereas the remaining measurements have weights $w = 2.5\text{e}^{3}$. We consider 24 operating conditions, which are generated by multiplying all the demands by the scale factors, and that the topology of the systems remains unchanged. Also, we set $\alpha = 0.3$ and $\beta = 1$ and $\beta = 1.5$ for the S-score and L-score, respectively.

\subsection{4-bus system}

\begin{figure}[t]
\begin{tikzpicture}[thick, scale=0.6,every label/.append style={font=\scriptsize}]
\newcommand{\gen}[3]{
  \node[circle,draw,thick,minimum width=5mm,inner sep=0pt,fill=white,#3] (#1) at (#2){};
  \node at ($(#2)-(0,0.5mm)$) {{\footnotesize $\sim$}} ;}
\newcommand{\Vm}[3]{
  \node[circle,draw,thick,minimum width=1.75mm,inner sep=0pt,fill=white,#3] (#1) at (#2){};}
\newcommand{\Pm}[3]{
 \node[mark size=3pt,color=black,#3] (#1) at (#2) {\pgfuseplotmark{triangle*}};}
\newcommand{\Qm}[3]{
 \node[mark size=3pt,color=gray,#3] (#1) at (#2) {\pgfuseplotmark{triangle*}};}  
\newcommand{\Pfm}[3]{
 \node[mark size=3pt,color=black,#3] (#1) at (#2) {\pgfuseplotmark{pentagon*}};}
\newcommand{\Qfm}[3]{
 \node[mark size=3pt,color=gray,#3] (#1) at (#2) {\pgfuseplotmark{pentagon*}};}  
\newcommand{\m}[3]{
 \node[mark size=0pt,color=white,#3] (#1) at (#2) {\pgfuseplotmark{pentagon*}};}   
  
\tikzstyle{line}=[-, thick]
\tikzstyle{loadline}=[->,thick,>=stealth']
\tikzstyle{busbar} = [rectangle,draw,fill=black,inner sep=0pt];
\tikzstyle{hbus} = [busbar,minimum width=5mm,minimum height=0.5pt];

\coordinate (c5) at (-2,0);
\coordinate (c6) at (1,0);
\coordinate (c7) at (6.5,0);

\coordinate (c4) at (4.5,2);

\coordinate (c1) at (0,6);
\coordinate (c2) at (4.5,6);
\coordinate (c3) at (9,6);

\node[hbus,minimum width=15mm,label=left:Bus 1] (b1) at (c1) {};
\node[hbus,minimum width=15mm,label=left:Bus 2] (b2) at (c2) {};
\node[hbus,minimum width=15mm,label=left:Bus 3] (b3) at (c3) {};
\node[hbus,minimum width=15mm,label=left:Bus 4] (b4) at (c4) {};

\draw[line] ([xshift=10mm] b1.south) -- ([xshift=10mm,yshift=-5mm] b1.south) -- ([xshift=-10mm, yshift=-5mm] b2.south) -- ([xshift=-10mm] b2.south); % Branch 1-2.
\draw[line] ([xshift=10mm] b2.south) |- ([xshift=-10mm, yshift=-5mm] b3.south) -- ([xshift=-10mm] b3.south); % Branch 2-5.
\draw[line] ([xshift=-10mm] b1.south) -- ([xshift=-10mm,yshift=-7.5mm] b1.south) -- ([xshift=-10mm, yshift=5mm] b4.north) -- ([xshift=-10mm] b4.north); % Branch 1-4.
\draw[line] ([xshift=10mm] b3.south) -- ([xshift=10mm,yshift=-7.5mm] b3.south) -- ([xshift=10mm, yshift=5mm] b4.north) -- ([xshift=10mm] b4.north); % Branch 1-4.

% Generator 1.
\gen{g1}{$(c1)+(0mm,12mm)$}{}
\draw[line] (b1.north) -- (g1.south);

\draw[loadline] (b3.north) -- ++(0,1) {};
\draw[loadline] (b4.south) -- ++(0,-1) {};

\Vm{V1}{$(c1)+(0mm,0mm)$}{label={[yshift=-0.65cm]$v_1^m$}}
\Vm{V2}{$(c2)+(0mm,0mm)$}{label={[yshift=-0.65cm]$v_2^m$}}
\Pm{P1}{$(c1)+(-4mm,3mm)$}{label={[xshift=-0.35cm, yshift=-0.3cm]$P_1^m$}}
\Pm{P2}{$(c2)+(-4mm,3mm)$}{label={[xshift=-0.35cm, yshift=-0.3cm]$P_2^m$}}
\Pm{P3}{$(c3)+(-4mm,3mm)$}{label={[xshift=-0.35cm, yshift=-0.3cm]$P_3^m$}}
\Qm{Q1}{$(c1)+(4mm,3mm)$}{label={[xshift=0.35cm, yshift=-0.3cm]$Q_1^m$}}
\Qm{Q2}{$(c2)+(4mm,3mm)$}{label={[xshift=0.35cm, yshift=-0.3cm]$Q_2^m$}}
\Qm{Q3}{$(c3)+(4mm,3mm)$}{label={[xshift=0.35cm, yshift=-0.3cm]$Q_3^m$}}
\m{P14}{$(c1)+(-14mm,-7mm)$}{label={[xshift=-0.3cm, yshift=-0.4cm]$P_{14}^m \pmb{\downarrow}$}}
\m{P32}{$(c3)+(-14mm,-9mm)$}{label={[xshift=-0.4cm, yshift=-0.4cm]$\pmb{\leftarrow} P_{32}^m$}}
\m{P34}{$(c3)+(6mm,-4mm)$}{label={[xshift=-0.3cm, yshift=-0.4cm]$P_{34}^m \pmb{\downarrow}$}}
\m{Q34}{$(c3)+(14mm,-4mm)$}{label={[xshift=0.3cm, yshift=-0.4cm]$\pmb{\Downarrow} Q_{34}^m $}}

\Vm{V}{$(c5)+(0mm,0mm)$}{label=right: Voltage}
\Pm{P}{$(c6)+(0mm,4mm)$}{label=right: Active power injection}
\Qm{Q}{$(c6)+(0mm,-4mm)$}{label=right: Reactive power injection}
\m{Pf}{$(c7)+(0mm,4mm)$}{label=right: $\pmb{\longrightarrow}$ Active power flow}
\m{Qf}{$(c7)+(0mm,-4mm)$}{label=right: $\pmb{\Longrightarrow}$  Reactive power flow}

\end{tikzpicture}
   \caption{4-bus system - One-line diagram and measurement location ($P_2^m$ and $Q_2^m$ are zero-injection measurements).}
   \label{fig:case4}
\end{figure}
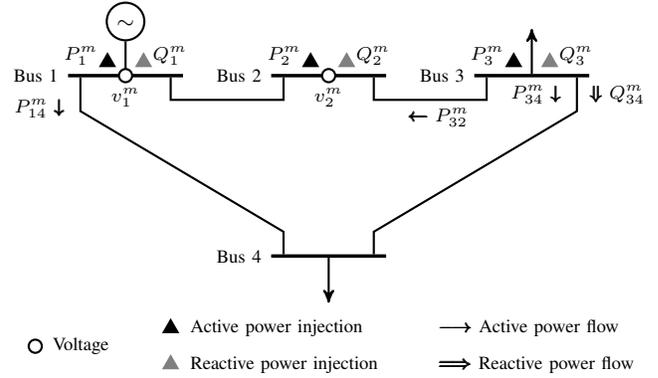

The 4-bus system and its measurement configuration are depicted in Fig.~\ref{fig:case4}. The presented measurement configuration provides a redundancy ratio of 1.71 where $P_2^m$ and $Q_2^m$ are zero-injection measurements.

\begin{figure}[t]
    \centering
    \includegraphics[width = 0.48\textwidth]{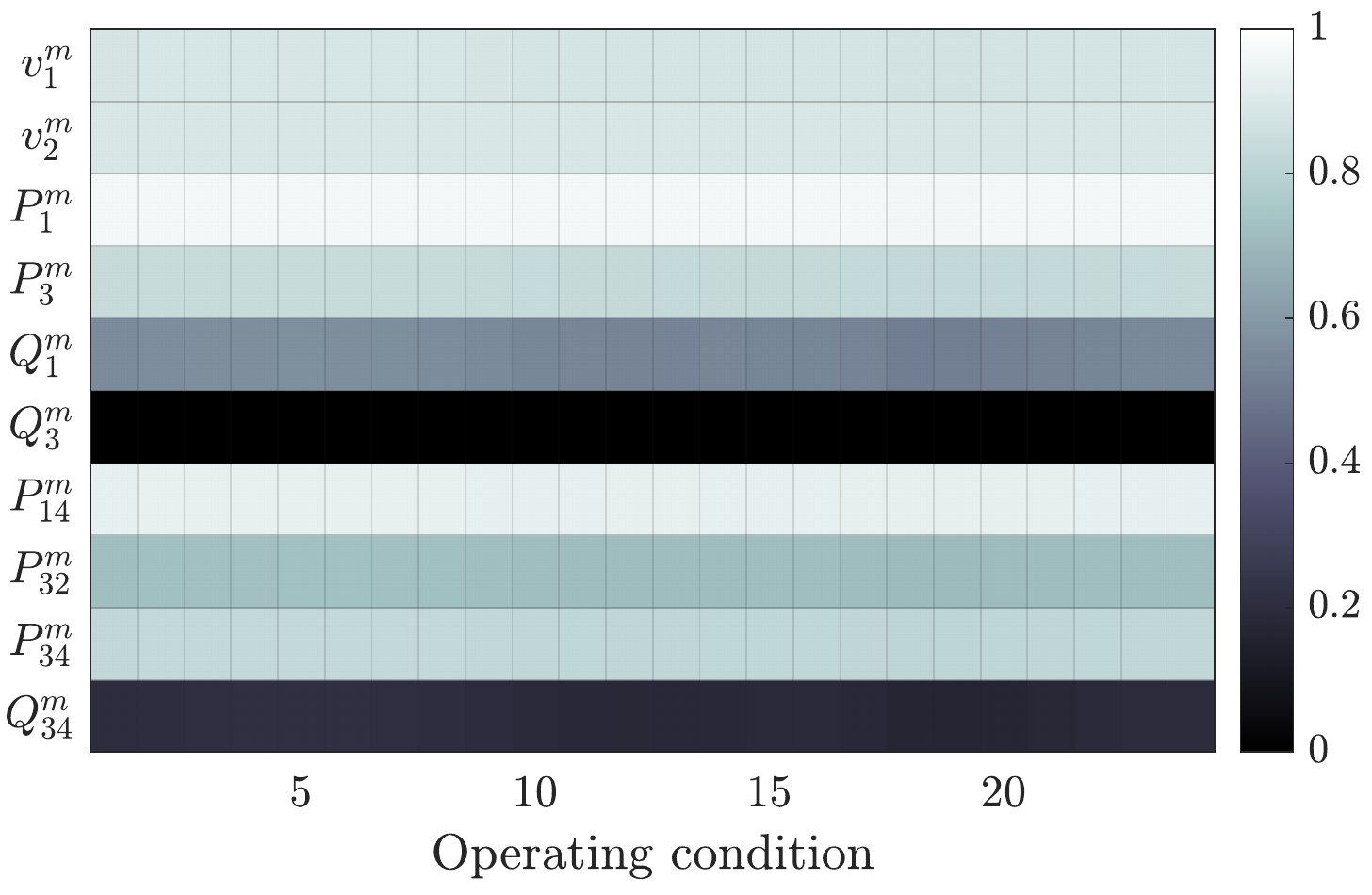}
    \caption{4-bus system - Normalized sensitivity of the objective with respect to measurements. The x-axis represents the 24 operating conditions. A darker color means that the measurement is more likely to be a target of undetected FDIAs.} 
    \label{fig:dJdA_4}
\end{figure}

The sensitivity of the objective function with respect to the measurements is depicted in Fig.~\ref{fig:dJdA_4}. We note that the magnitude of the sensitivities remains almost invariant with the operating conditions. $P_{1}^m$ is the measurement with the largest normalized sensitivity $\left( \left\lvert z_k \frac{\partial J}{\partial z_{\ell}}\right\rvert  / \max\limits_k \left\{ \left\lvert z_k \frac{\partial J}{\partial z_k}\right\rvert \right\} \right)$, thus the less vulnerable in terms of stealthiness as an FDIA against it is unlikely to remain undetected. Conversely, $Q_3^m$ is the measurement with the smallest sensitivity followed by $Q_{34}^m$, $Q_{1}^m$ and $P_{32}^m$, respectively. 

\begin{figure}[t]
    \centering
    \includegraphics[width = 0.48\textwidth]{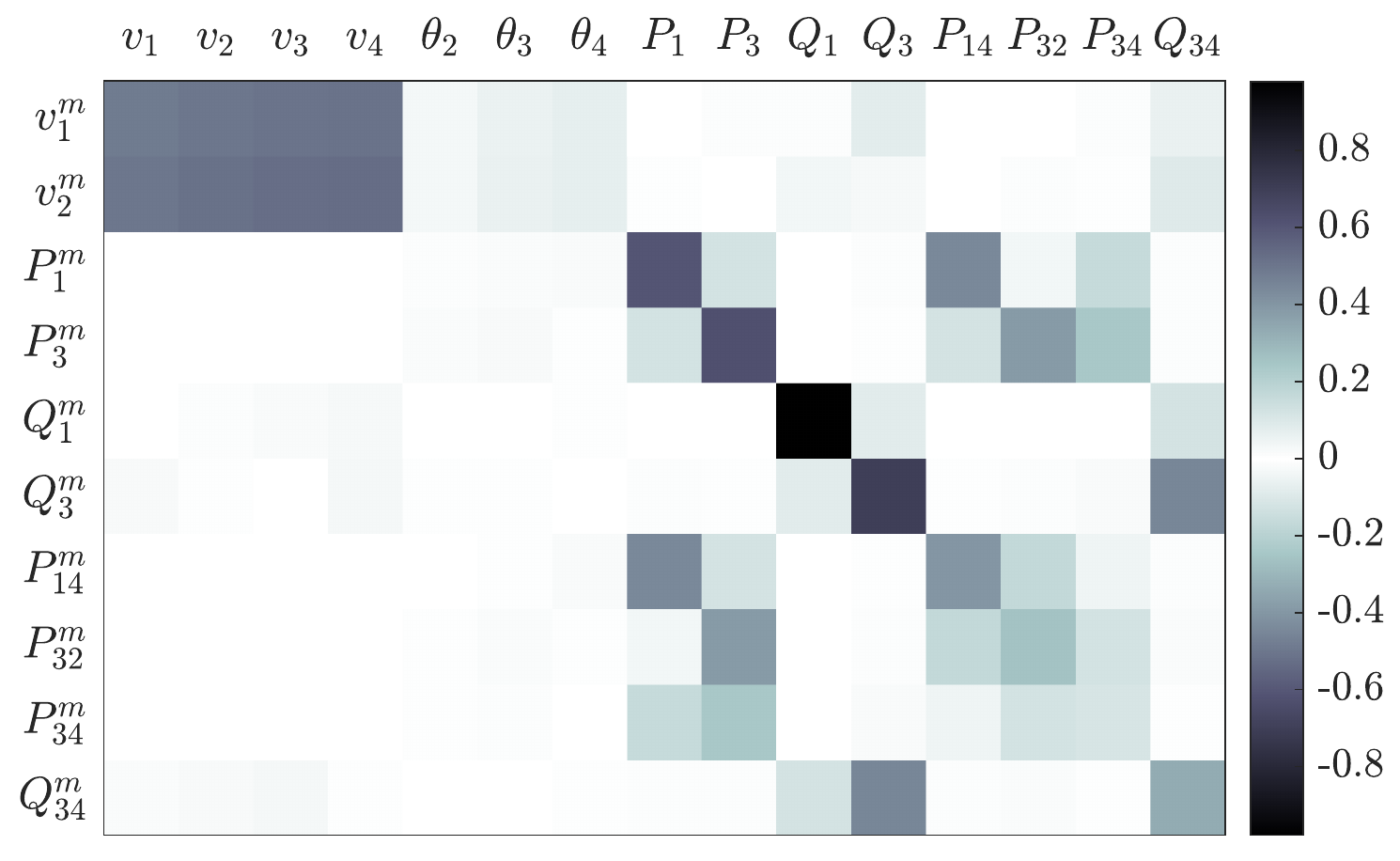}
    \caption{4-bus system - Sensitivity of estimated variables with respect to measurements (Scale factor = 1). The sensitivities with larger absolute value are depicted with darker colors. A darker color means that the measurement is more likely to be a target of impactful FDIAs.}
    \label{fig:dXdA_4}
\end{figure}

Likewise, the sensitivity of the variable estimates with respect to the measurements when at maximum demand is depicted in Fig.~\ref{fig:dXdA_4}. The measurements with the highest self-sensitivities are $Q_1^m,\,Q_3^m,\,P_3^m$ and $P_1^m$, respectively, which are the most vulnerable ones in terms of impactfulness. Specifically, $Q_1^m$ shows the largest self-sensitivity. $v_1^m$ and $v_2^m$ show the largest impact on the other variable estimates ($v_3$ and $v_4$).
A FDIA compromising these measurements will have a significant impact on the corresponding variable estimates. Furthermore, it is convenient to analyze the dependence of the variable estimates with respect to each measurement. $v_1^m$ and $v_2^m$ show important influence on the estimates of $v_3$ and $v_4$. Similarly, $P_1^m$ and $Q_3^m$ (resp. $P_{14}^m$ and $Q_{34}^m$  have a significant influence on the variable estimates of $P_{14}$ and $Q_{34}$ (resp. $P_1$ and $Q_3$), respectively.

Table~\ref{tab:sensit_4} provides the three proposed scores of each measurement at a scale factor equal to 1. $Q_1^m$ is the most vulnerable measurement followed by $Q_3^m$. All these measurements are vulnerable  due to the lack of redundancy of reactive power measurements near buses 1 and 3. These measurements exhibit the highest potential to be targeted for cyber-attacks. The best chance of staging a stealthy and impactful FDIA if corrupting any of these two measurements.

\begin{figure}[t]
    \centering
    \includegraphics[width = 0.49\textwidth]{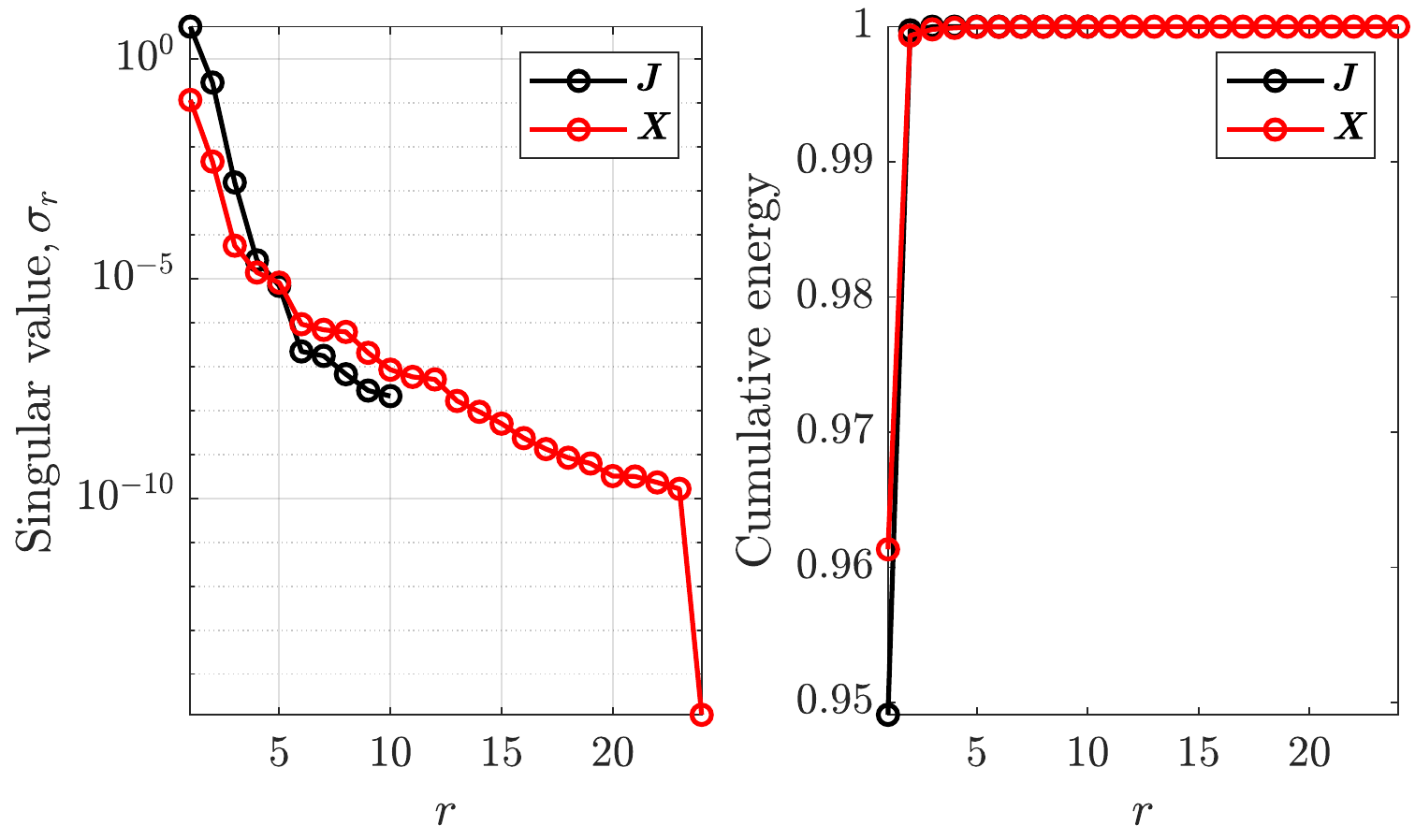}
    \caption{4-bus system - (Left) Singular values $\sigma_r$. (Right) Cumulative energy in the first $r$ singular values.}
    \label{fig:SVD_4}
\end{figure}

Fig.~\ref{fig:SVD_4} shows that the leading singular value of $\tilde{\boldsymbol{X}}$ and $\tilde{\boldsymbol{J}}$ account for almost $95\%$ and $96\%$ of their variance. Thus, both sensitivity vectors are almost invariant to the different operating points, which lead us to the conclusion that the vulnerabilities are mainly dependent on the network topology and its parameters, and the measurements' configuration.

\begin{table}[t]
\caption{4-bus system - Vulnerability Scores}
\renewcommand{\arraystretch}{1.4}
\centering
\label{tab:sensit_4}
%\resizebox{\columnwidth}{!}{
\begin{tabular}{c c c  c}
\hline 
\hline
Meas. & S-score & L-score & V-score \\
\hline
$v_1^m$    & 0.4037 & 1.0000 & 0.8211 \\
$v_2^m$    & 0.3506 & 1.0000 & 0.8052 \\
$P_{1}^m$  & 0.0334 & 0.8692 & 0.6184 \\
$P_{3}^m$  & 0.5092 & 0.8794 & 0.7683 \\
$Q_{1}^m$  & 0.7913 & 0.9989 & 0.9366 \\
$Q_{3}^m$  & 0.8927 & 0.9188 & 0.9109 \\
$P_{14}^m$ & 0.1974 & 0.6919 & 0.5436 \\
$P_{32}^m$ & 0.6747 & 0.5112 & 0.5603 \\
$P_{34}^m$ & 0.5331 & 0.2678 & 0.3474 \\
$Q_{34}^m$ & 0.8725 & 0.6076 & 0.6871 \\
\hline 
\hline
\end{tabular}
\end{table}

\subsection{New England 39-bus system}
We consider that the New England 39-bus system has the following measurements: all the voltage magnitudes, active and reactive power injections at all generating buses, and the active and reactive power flows at the sending ends of all lines. The system data can be retrieved from MATPOWER \cite{Zimmerman2011}.

\begin{figure}[t]
    \centering
    \includegraphics[width = 0.48\textwidth]{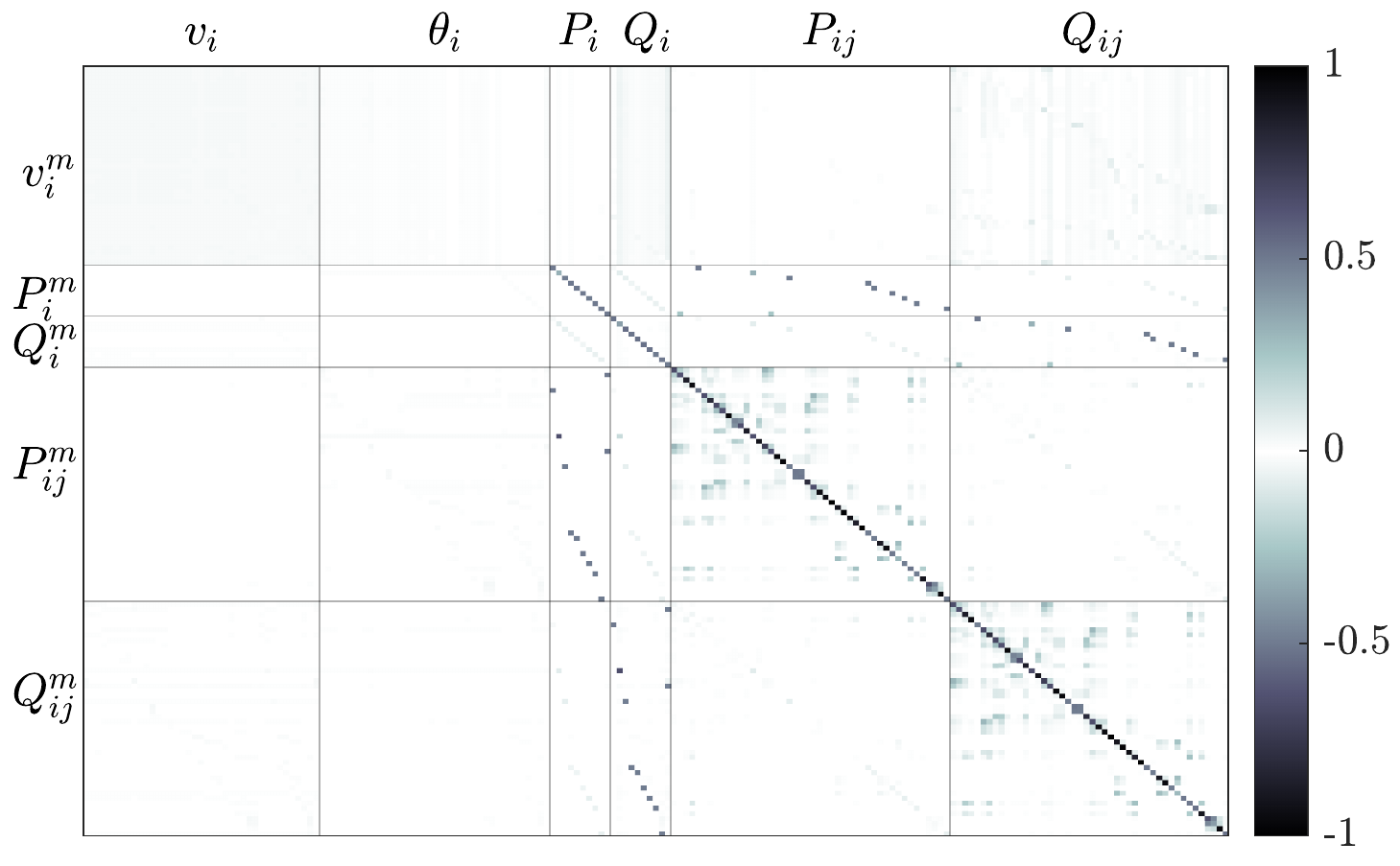}
    \caption{New England 39-bus system - Sensitivity of estimated variables with respect to measurements (Scale factor = 1).}
    \label{fig:dXdA_39}
\end{figure}

The sensitivities of the estimated variables with respect to measurements are depicted in Fig. \ref{fig:dXdA_39}. The voltage measurements are not leverage points since their self-sensitivities are small. Conversely, the majority of the active and reactive power flow and injection variables have high sensitivity with respect to their corresponding measurements.

\begin{figure*}[t]
    \centering
    \includegraphics[width = \textwidth, trim={2.5cm 0 2.5cm 0},clip]{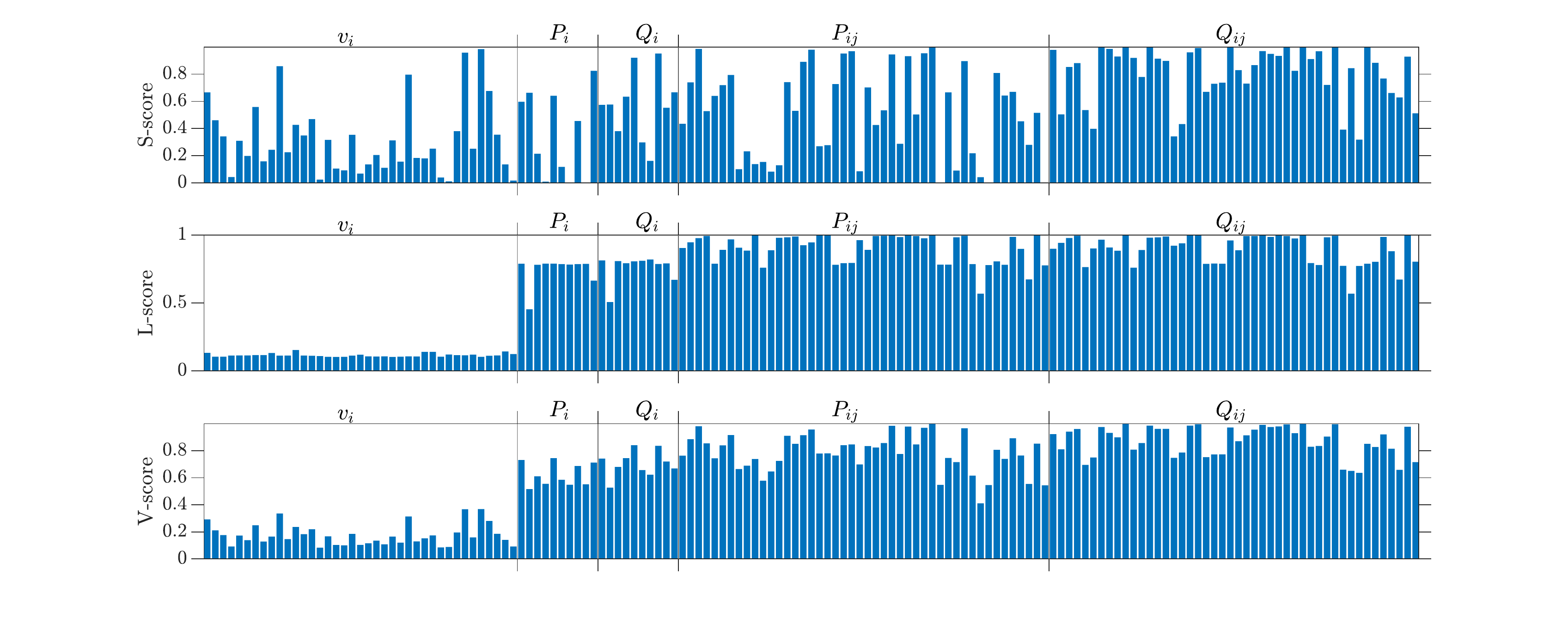}
    \vspace{-1cm}
    \caption{New England 39-bus system - Vulnerability scores.}
    \label{fig:Scores39}
\end{figure*}

Additionally, Fig. \ref{fig:Scores39} depicts the vulnerability scores of all the measurements. An important number of measurements show an S-score near 1, which makes them attractive to the adversary from the stealthiness point of view. The L-score does not show the same distribution; however, there are 48 measurements whose L-score is greater than 0.8. Finally, 22 active and reactive power flow measurements have a V-score greater than 0.95. These results show that the lack of redundancy of active and reactive power measurements is not localized in a certain area of the system, which may be due to the low redundancy ratio.

\begin{figure}
    \centering
    \includegraphics[width = 0.49\textwidth, trim={2.75cm 9.15cm 2.75cm 9.15cm},clip]{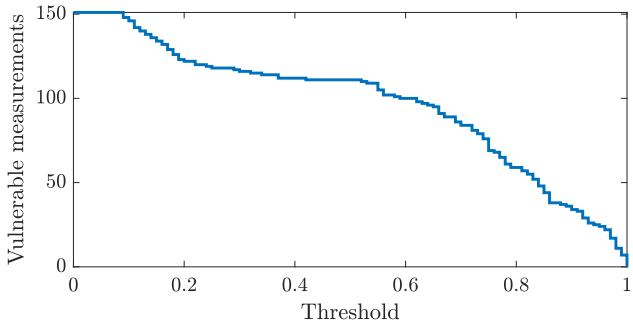}
    \caption{New England 39-bus system - Vulnerable measurements. A measurement $z_\ell$ is considered vulnerable if its $\text{V-score}(z_\ell) \ge \text{threshold}$.}
    \label{fig:vuln_39}
\end{figure}

We also provide the number of vulnerable measurements as a function of different threshold values in Fig.~\ref{fig:vuln_39}. We consider that a measurement $z_\ell$ is vulnerable if its $\text{V-score}(z_\ell) \ge \text{threshold}$. A smaller threshold implies higher conservativeness as it results in declaring a larger number of measurements as vulnerable. Table \ref{tab:vuln_39} lists the ten most critical measurements, $\text{V-score}\ge0.9836$, in descending order of their V-score. Clearly, these measurements are potential targets of FDIAs as their scores are close to 1, which means that if they are perturbed, they significantly influence their corresponding variable estimates.

\begin{table}[t]
\caption{39-bus system - Vulnerability Scores}
\renewcommand{\arraystretch}{1.45}
\centering
\label{tab:vuln_39}
%\resizebox{\columnwidth}{!}{
\begin{tabular}{c c c c }
\hline 
\hline
Meas. & S-score & L-score & V-score \\
\hline
$P_{19-20}^m$ & 0.999176 & 1        & 0.999753 \\
$Q_{19-20}^m$ & 0.997233 & 0.999953 & 0.999137 \\
$Q_{5-6}^m$   & 0.998304 & 0.996979 & 0.997377 \\
$Q_{22-23}^m$ & 0.996568 & 0.995312 & 0.995689 \\
$Q_{10-13}^m$ & 0.992317 & 0.997078 & 0.995649 \\
$Q_{17-18}^m$ & 0.996833 & 0.992793 & 0.994005 \\
$Q_{16-19}^m$ & 0.969506 & 0.999661 & 0.990615 \\
$Q_{6-11}^m$  & 0.99855  & 0.980853 & 0.986162 \\
$Q_{10-11}^m$ & 0.960126 & 0.997099 & 0.986007 \\
$P_{16-19}^m$ & 0.945391 & 1        & 0.983617 \\
\hline 
\hline
\end{tabular}
\end{table}

The leading singular values of $\tilde{\boldsymbol{X}}$ and $\tilde{\boldsymbol{J}}$, presented in Table \ref{tab:SVD_39}, capture around $65\%$ and $80\%$ of their variance, respectively. In the case of $\tilde{\boldsymbol{J}}$, its first three singular values account for more than $97\%$ of its variance. On the other hand, the four leading singular values of $\tilde{\boldsymbol{X}}$ capture around $90\%$ of its variance. These results show the low-rank characteristic of both matrices.

\begin{table}[t]
\caption{New England 39-bus system - SVD results}
\renewcommand{\arraystretch}{1.4}
\centering
\label{tab:SVD_39}
%\resizebox{\columnwidth}{!}{
\begin{tabular}{c|c c || c c}
\hline 
\hline
\multirow{2}{*}{$r$}& \multicolumn{2}{c||}{Singular values, $\sigma_r$} & \multicolumn{2}{c}{Cumulative energy}  \\
 & $\boldsymbol{J}$ & $\boldsymbol{X}$ & $\boldsymbol{J}$ & $\boldsymbol{X}$\\
\hline
1 & 304.2029 & 3.4085 & 0.8010 & 0.6525\\
2 & 54.4586 & 1.2666 &  0.9443 & 0.8950\\
3 & 14.2156  & 0.4722 & 0.9818 & 0.9854\\
4 & 4.6080  & 0.0643 &  0.9939 & 0.9977\\
5 & 2.2604  & 0.0105 &  0.9999 & 0.9997\\
6 & 0.0398  & 0.0011 &  1.000 & 0.9999\\
\hline 
\hline
\end{tabular}%}
\end{table}

\section{Conclusions}
This paper proposes a technique based on sensitivity analysis to identify measurements with a high potential of being the target of false data injection attacks. We characterize the vulnerability of each measurement as a function of their potential to impact the variable estimates and to remain stealthy.

In our numerical studies, we demonstrate that there is a subset of measurements that show both characteristics, thus being the most vulnerable to FDIAs.  Furthermore, we numerically demonstrate that such vulnerabilities remain almost invariant to the system's operating condition, which implies that they are mainly dependent on the network topology and its parameters, and the measurements' configuration. 

The proposed technique can be used to identify the most vulnerable measurements. Additionally, identifying such measurements can be used to secure their telemetry or to allocate new measurements to improve the redundancy over specific areas of the system. %These countermeasures can increase the system security.

\section*{Appendix A} \label{app:A}
%\subsection{Appendix A} \label{app:A}
The auxiliary submatrices and vectors in \eqref{eq:8} necessary for computing the sensitivities are defined below:
\begin{gather}
\boldsymbol{J}_{\boldsymbol{x}(1\times n)} = \left[ \nabla_{\boldsymbol{x}}J(\boldsymbol{x}^\ast,\boldsymbol{a},\boldsymbol{z})\right]^\top, \\
\boldsymbol{J}_{\boldsymbol{a}(1\times q)} = \left[ \nabla_{\boldsymbol{a}}J(\boldsymbol{x}^\ast,\boldsymbol{a},\boldsymbol{z})\right]^\top, \\
\boldsymbol{J}_{\boldsymbol{z}(1\times p)} = \left[ \nabla_{\boldsymbol{z}}J(\boldsymbol{x}^\ast,\boldsymbol{a},\boldsymbol{z})\right]^\top, \\
\boldsymbol{J}_{\boldsymbol{xx}(n\times n)} =\nabla_{\boldsymbol{xx}}J(\boldsymbol{x}^\ast,\boldsymbol{a},\boldsymbol{z}) + \sum_{k=1}^{r}\lambda_k^\ast \nabla_{\boldsymbol{xx}} c_k(\boldsymbol{x}^\ast,\boldsymbol{a}), \\
\boldsymbol{J}_{\boldsymbol{xa}(n\times q)} = \nabla_{\boldsymbol{xa}}J(\boldsymbol{x}^\ast,\boldsymbol{a},\boldsymbol{z}) + \sum_{k=1}^{r}\lambda_k^\ast \nabla_{\boldsymbol{xa}} c_k(\boldsymbol{x}^\ast,\boldsymbol{a}), \\
\boldsymbol{J}_{\boldsymbol{xz}(n\times p)} = \nabla_{\boldsymbol{xz}}J(\boldsymbol{x}^\ast,\boldsymbol{a},\boldsymbol{z}) d\boldsymbol{z}, \\
\boldsymbol{C}_{\boldsymbol{x}(r\times n)} =\left[\nabla_{\boldsymbol{x}} \boldsymbol{c}(\boldsymbol{x}^\ast,\boldsymbol{a})\right]^\top, \\ 
\boldsymbol{C}_{\boldsymbol{a}(r\times q)} = \left[\nabla_{\boldsymbol{a}} \boldsymbol{c}(\boldsymbol{x}^\ast,\boldsymbol{a})\right]^\top.
\end{gather}

\begin{IEEEbiography}[{\includegraphics[width=1in,height=1.25in,clip,keepaspectratio]{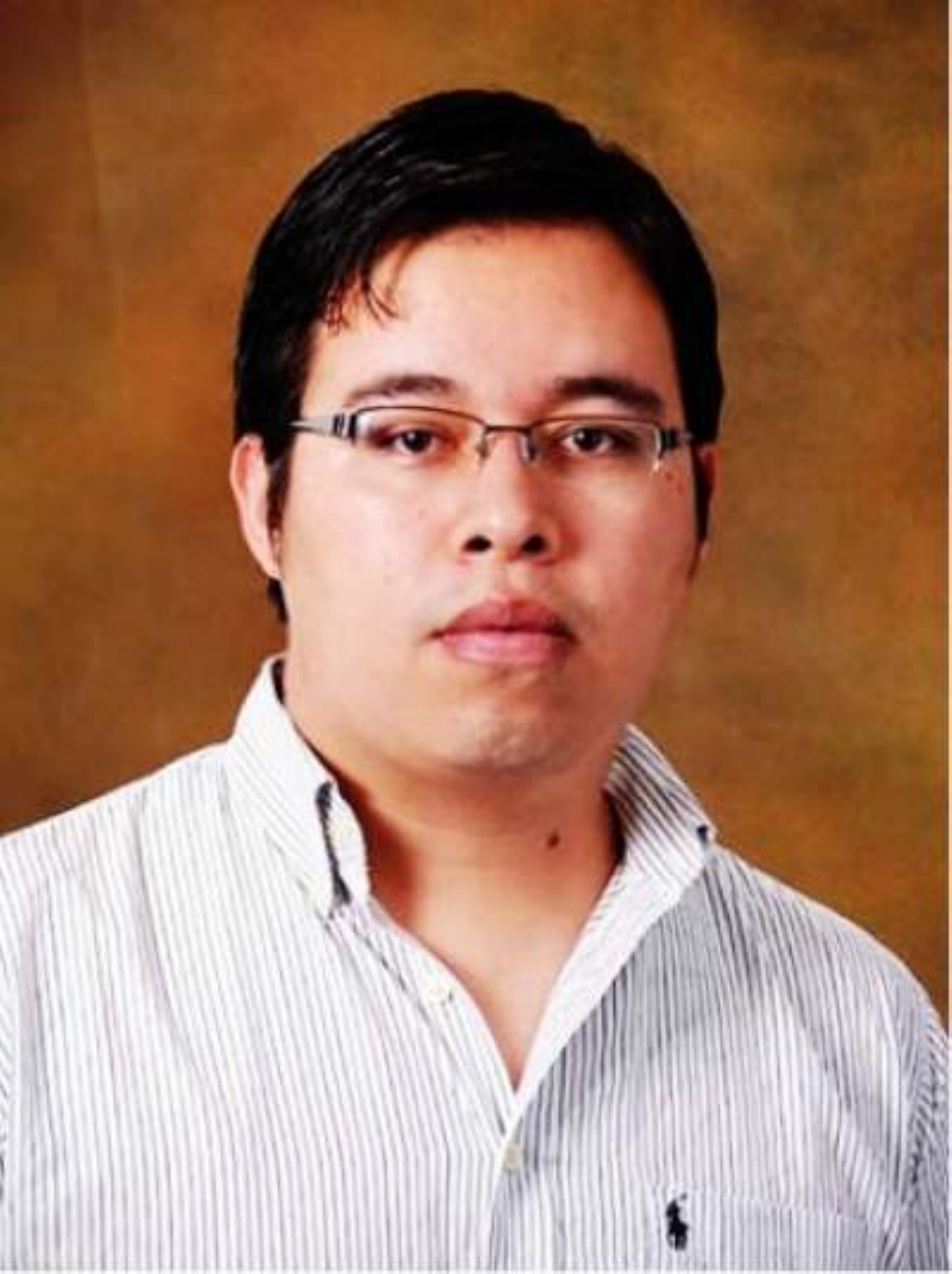}}]{Gonzalo E. Constante-Flores}
    (S'$13$) received the electrical engineering degree from Escuela Polit\'ecnica Nacional (EPN), Quito, Ecuador, in 2014, and the M.Sc. degree from The Ohio State University, Columbus, OH, in 2018, where he is currently working towards the Ph.D. degree in the Department of Electrical and Computer Engineering.\\\indent From 2013 to 2016, he was with the Department of Electrical Energy at EPN. His research interests include optimization and control of power systems, electricity markets, and the integration of renewable energy and electric vehicles into electric power systems.
\end{IEEEbiography}
\begin{IEEEbiography}
	[{\includegraphics[clip,width=1in,height=1.25in,keepaspectratio]{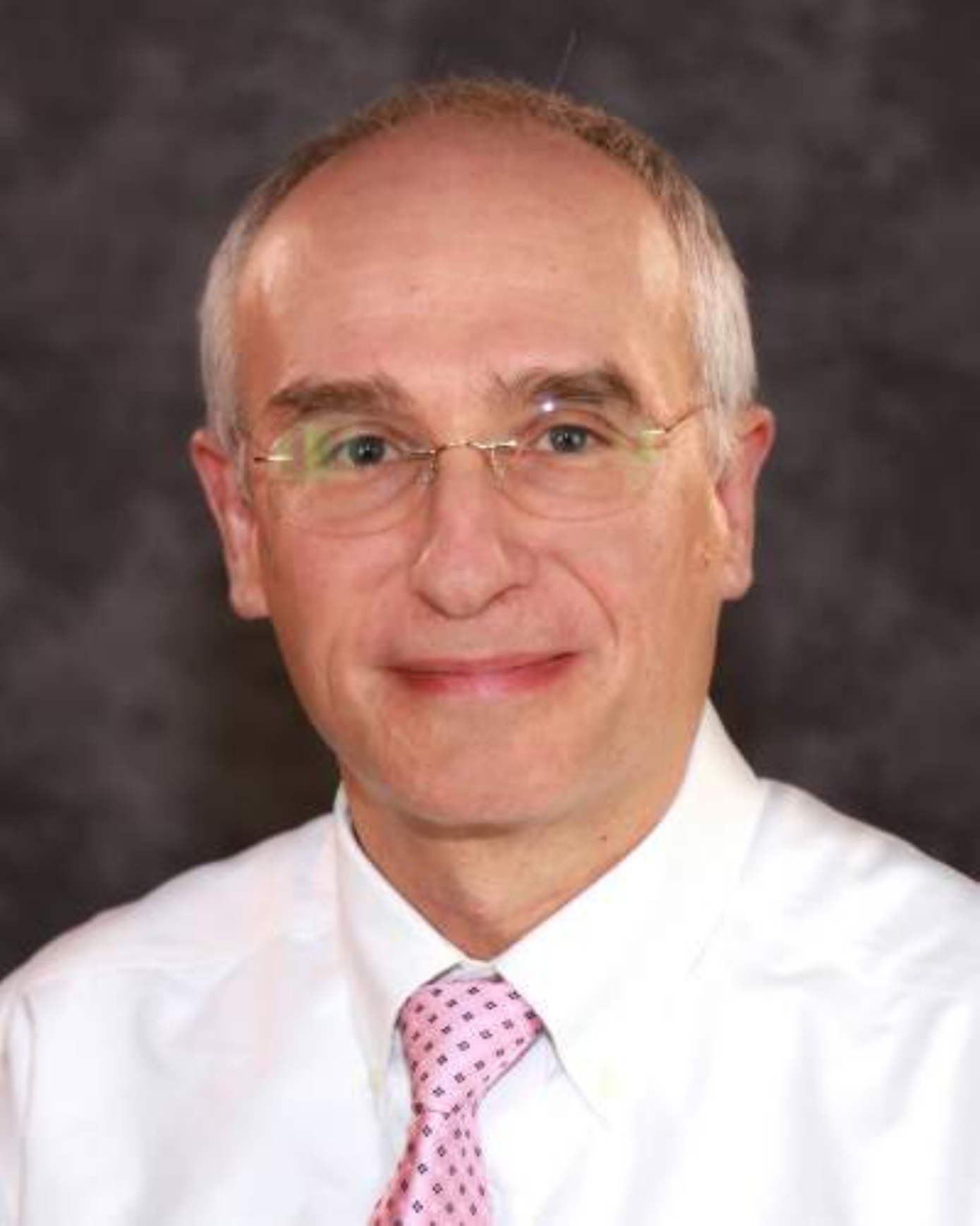}}]
	{Antonio J. Conejo}
	(F'$04$) received the M.S.\ degree from the Massachusetts Institute of Technology, Cambridge, MA, in~$1987$, and the Ph.D.\ degree from the Royal Institute of Technology, Stockholm, Sweden, in~$1990$.\\\indent He is currently a professor in the Department of Integrated Systems Engineering and the Department of Electrical and Computer Engineering, The Ohio State University, Columbus, OH. His research interests include control, operations, planning, economics and regulation of electric energy systems, as well as statistics and optimization theory and its applications.
\end{IEEEbiography}
\begin{IEEEbiography}
	[{\includegraphics[clip,width=1in,height=1.25in,keepaspectratio]{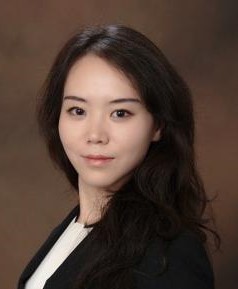}}]{Jiankang Wang}
	(M'$09$) received the M.S. and Ph.D. degree from MIT, Cambridge, MA in 2009 and 2013, respectively. \\\indent She joined the Department of Electrical and Computer Engineering at The Ohio State University, as an assistant professor in 2014, and appointed as the lead technical specialist of California ISO in 2018. Her research areas include power system cybersecurity, renewable energy integration, and electricity markets.  
\end{IEEEbiography}

\end{document}